\documentclass[11pt]{article}

\usepackage{amssymb,amsmath,latexsym}
\usepackage[dvips]{graphics}

\def\EE{{\cal E}}
\def\FF{{\cal F}}

\def\CC{{\cal C}}
\def\DD{{\cal D}}
\def\LL{{\cal L}}

\def\P{{\bf P}}

\def\RR{\mathbb{R}}
\def\N{\mathbb{N}}
\def\HH{\mathbb{H}}

\def\p{{\bf p}}
\def\n{{\bf n}}
\def\D{{\bf D}}
\def\rf{{\rm ref}}

\def\<{\langle}
\def\>{\rangle}
\def\pf{\noindent{\bf Proof.} }
\def\eps{{\varepsilon}}

\def\wt{\widetilde}
\def\wh{\widehat}
\def\dis{\displaystyle}
\def\bl{{\rm BL}}

\def\qed{{\hfill $\Box$ \bigskip}}

\numberwithin{equation}{section}

\newtheorem{thm}{Theorem}[section]
\newtheorem{prop}[thm]{Proposition}
\newtheorem{remark}[thm]{Remark}
\newtheorem{defn}[thm]{Definition}
\newtheorem{corollary}[thm]{Corollary}

\newtheorem{lemma}[thm]{Lemma}

\begin{document}

\title{\bf On unique extension of time changed
reflecting Brownian motions \medskip \\
\rm \small (To appear in {\it Ann. Inst. Henri Poincare Probab.
Statist.})}

\bigskip

\author{{\bf Zhen-Qing Chen}\thanks{Research partially supported
by NSF Grant  DMS-06000206.} \  \quad  and \quad {\bf Masatoshi
Fukushima}\thanks{Research supported by Grant-in-Aid for Scientific
Research of MEXT No.19540125}}

\date{(August 26,  2008)}

\maketitle

\begin{abstract} Let $D$ be an unbounded domain in $\RR^d$ with $d\geq 3$. We show
 that if $D$ contains an unbounded uniform domain, then
 the symmetric reflecting Brownian motion (RBM)  on $\overline D$
is transient. Next assume that RBM $X$ on $\overline D$ is transient and
 let $Y$ be its time change by Revuz measure ${\bf 1}_D(x) m(x)dx$
 for a strictly positive continuous integrable function $m$ on $\overline D$.
 We further show that if there is some $r>0$ so that $D\setminus \overline {B(0, r)}$ is
an unbounded uniform domain, then $Y$ admits one and only one symmetric diffusion
that genuinely extends it
 and admits no killings.
 In other words, in this case $X$ (or equivalently, $Y$)
has a unique Martin boundary point at infinity.

\bigskip

 \centerline{\bf R\'esum\'e}

 \medskip

Notons $D$ une domaine non born{\'e}e dans $\mathbb R^d$ avec $d \ge 3$.
Nous montrons que si $D$  contient une domaine uniforme non born\'ee,
alors le mouvement brownien refl{\'e}tent (RBM) sur $\overline{D}$ est
transitoire.  Suivant nous supposons que RBM $X$ sur $\overline{D}$ est
transitoire et notons $Y$ sa changement de temps par une mesure Revuz
$\textbf{1}_D(x)m(x)dx$ pour  un fonction $m$ strictement positif continue
int{\'e}grable sur $\overline{D}$.  En plus montrons que s'il existe un $r >
0$ telles que $D \setminus \overline{B(0,r)}$ est une domaine uniforme non born\'ee,
alors $Y$ admet un et seulement un diffusion symmetric que
l'{\'e}tend et admet pas les meurtres.  Autrement dit, dans ce cas $X$ (ou
{\'e}galement $Y$) a un point de bord Martin unicit{\'e} {\`a} l'infini.
\end{abstract}

\vspace{.3truein}


\bigskip


\section{Introduction}

Consider an unbounded domain $D$ in $\RR^d$ with $d\ge 3$ possessing continuous boundary.
Then $(\frac12\D, W^{1,2}(D))$ is a regular strongly local Dirichlet form on $L^2(\overline D; dx)$, where $\D(u,v)= \int_D \nabla u(x) \cdot \nabla v(x) \, dx$. The symmetric diffusion process $X$ associated with this form
is the symmetric reflecting Brownian motion on $\overline D$ with infinite lifetime.
In this paper, we are concerned with the following two questions.

\begin{description}
\item{(i)} When $X$ is transient?

\item{(ii)} If $X$ is transient,
 then almost all sample paths of $X$ approach to the point at infinity $\partial$ of $\overline{D}$ as time goes to infinity.
 Let $Y$ be its time change by Revuz measure ${\bf 1}_D(x) m(x)dx$
 for a strictly positive continuous integrable function $m$ on $\overline D$.
  The time-changed process $Y$ has finite lifetime with positive probability.
 How many symmetric diffusions are there that genuinely extends $Y$?
\end{description}

 We show that  if  $D$ contains
an unbounded uniform domain (see Definition \ref{D:2.1} below), then RBM $X$ on $\overline D$
is transient. We further show that if there is some $r>0$ so that $D\setminus \overline {B(0, r)}$
is an unbounded uniform domain, then $Y$ admits a unique genuine diffusion extension
that admits no killings.
Here $B(x, r)$ denotes the ball in $\RR^d$ centered at $x$ with radius $r$.
 The extension diffusion can be obtained through one-point darning of $Y$ at $\partial.$
  It can also be obtained through the active reflected Dirichlet form of $Y$.
  The key of
  our approach is to identify the
 reflected Dirichlet space of $Y$
   with
 the space
 $$ \bl(D) :=
 \Big\{u\in L^2_{loc}(D):\frac{\partial u}{\partial x_i}\in L^2(D),\ 1\le i\le d \Big\}
 $$
  of
  Beppo Levi functions on $D$ and to show that $\bl (D)$ is the linear space spanned by $W_e^{1,2}(D) $
  and constant functions
  under suitable conditions on $D$.
Here $W^{1,2}_e (D)$ is the extended Dirichlet space of $(\tfrac12 \D, \, W^{1,2}(D))$;
see next section for its definition.

\section{Transience of reflecting Brownian motion}
Let $E$ be a locally compact separable metric space and
 $m$
a positive Radon measure on $E$ with full support.
 Numerical functions $f,\;g$ on $E$ is said to be $m$-equivalent if $m(f\neq g)=0$ we write as $f=g\ [m]$ in this case.
Let $(\EE,\FF)$ be a Dirichlet form on
 $L^2(E; m)$
and $(\FF_e,\EE)$ be its extended Dirichlet space.
 A function $u$ is in $\FF_e$
if and only if $|u|<\infty\ [m]$ and there exists $\{u_n\}\subset
\FF$ called an approximating sequence of $u$ such that $\{u_n\}$ is
$\EE$-Cauchy and $\lim_{n\to\infty}u_n=u\ [m].$  It holds then that
$\EE(u,u)=\lim_{n\to\infty}\EE(u_n,u_n).$  We know that
$\FF=\FF_e\cap L^2(E;m)$ and that every normal contraction operates
on $(\FF_e,\EE).$

Denote by $\{T_t;t>0\}$ be the $L^2$-semigroup generated by the Dirichlet form $\EE$ and
 define $\{ S_t,\; t>0\}$ by the Bochner integral $S_tf=\int_0^t T_sfds,\ f\in L^2(E;m).$
 The operator $S_t$
 then extends from $L^2\cap L^1$ to a bounded linear operator on $L^1(E;m)$
 such that $0\le S_tf\le S_{t'}f
 $ for $ 0\leq t<t'$ and $\ f\in L^1_+(E;m)$.
  Hence for $f\in L^1_+(E;m)$,
  \[ Gf(x)=\lim_{N\to\infty}S_Nf(x)( \le \infty) \]
  defines
  a function $Gf $
 uniquely up to the $m$-equivalence.
 The Dirichlet form $(\EE, \FF)$
is called {\it transient} (resp. {\it recurrent}) if $Gf<\infty\
[m]$ for some $f\in L_+^1(E;m)$ with $f>0\ [m]$ (resp.
$m(0<Gf<\infty)=0
 $ for every $ f\in L_+^1(E;m)$).
 An $m$-measurable
set $A\subset E$ is said to be $\{T_t\}$-{\it invariant} if
$1_AT_t1_{A^c}=0\ [m].$
 The Dirichlet form $(\EE,  \FF)$
is called
 $m$-{\it irreducible}
 if
 for any $\{T_t\}$-invariant set $A$, we have either $m(A)=0$ or $m(A^c)=0$.
 If the Dirichlet form $(\EE, \FF)$
is irreducible, then it is either transient or recurrent. The
following criteria are known (\cite{FOT}).

\noindent The Dirichlet form $(\EE, \FF)$ is transient if and only
if one of the next conditions is satisfied:
\begin{equation}\label{eqn:2.1}
(\FF_e,\EE)\ \hbox{\rm is a real Hilbert space};
\end{equation}
\begin{equation}\label{eqn:2.2}
u\in \FF_e,\quad \EE(u,u)=0\ \Longrightarrow\ u=0\ [m].
\end{equation}
The Dirichlet form $(\EE, \FF)$ is recurrent if and only if
\begin{equation}\label{eqn:2.3}
1\in \FF_e \quad  \hbox{ and } \quad \EE(1,1)=0.
\end{equation}
Denote by
 $C_c(E)$
the space of all continuous functions on $E$ with compact support.
  For a Dirichlet form $(\EE,\FF),$ a subspace $\CC$ of $\FF\cap C_c(E)$ is called a
  {\it core} of $\EE$ if $\CC$ is $\EE_1$-dense in $\FF$ and uniformly dense in $C_c(E).$
  A Dirichlet form $\EE$ is called {\it regular} if it possesses a core.
  A Dirichlet form $(\EE,\FF)$ is said to be {\it strongly local} if
$\EE(u,v)=0$ whenever $K={\rm Supp}[u\cdot m]$ is compact and $v$ is
constant in a neighborhood of $K.$ Any strongly local regular
Dirichlet form on $L^2(E;m)$ is known to admit an associated
$m$-symmetric diffusion process on
$E$ that admits no killing inside $E$.

In the following, we will take the state space $E$ to be a Euclidean domain $D$
in $\RR^d$.
 We denote by $L^p(D),\ p\ge 1,$ the $L^p$-space
of functions on $D$
 with respect to
the Lebesgue measure $dx$.  We focus our attention on the space
\begin{equation}\label{eqn:2.4}
\bl(D)=\left\{T: \frac{\partial T}{\partial x_i}\in L^2(D),\ 1\le i\le d\right\}\end{equation}
of Schwartz distributions $T.$  It is known
that any distribution $T\in \bl(D)$ can be identified with a function in $L_{loc}^2(D)$ (cf. L. Schwartz \cite{Schw}, J. Deny and J.L. Lions \cite{DL}) so that
\begin{equation}\label{eqn:e.2.5}
\bl(D)=\left\{u\in L^2_{loc}(D):\frac{\partial u}{\partial x_i}\in
L^2(D),\ 1\le i\le d \right\},
\end{equation}
where the derivatives are taken in Schwartz distribution sense.  Members in $\bl
(D)$ are called BL (Beppo Levi) {\it functions} on $D$.   For $u,v\in \bl(D),$ we put
\begin{equation}\label{eqn:2.6}
{\bf D}(u,v)=\sum_{i=1}^d\int_D \frac{\partial u}{\partial
x_i}\frac{\partial v}{\partial x_i}dx.
\end{equation}

The space $\bl(D)$ is known to enjoy the following properties (cf. \cite{DL}):
\begin{description}
\item[$(\bl.1)$] The quotient space $\dot{\bl}(D)$ of $\bl(D)$ by the subspace of constant functions is a Hilbert space with inner product ${\bf D}.$  Any ${\bf D}$-Cauchy sequence $u_n\in \bl(D)$ admits $u\in \bl(D)$ and constants $c_n$ such that $u_n$ is ${\bf D}$-convergent to $u$ and $u_n+c_n$ is $L^2_{loc}$-convergent to $u.$
\item[$(\bl.2)$] A function $u$ on $D$ is in $\bl(D)$ if and only if, for each $i\ (1\le i\le d),$ there is a version $u^{(i)}$ of $u$ such that it is absolutely continuous on almost all straight lines parallel to $x_i$-axis and the derivative $\partial u^{(i)}/\partial x_i$ in the ordinary sense (which exists a.e. on $D$) is in $L^2(D)$.  In this case, the ordinary derivatives coincide with the distribution derivatives of $u.$
\end{description}

The {\it Sobolev space of order}
 $(1, 2)$
on the domain $D\subset \RR^d$ is defined by
\begin{equation}\label{eqn:2.7}
 W^{1,2}(D)
=\bl(D)\cap L^2(D).
\end{equation}
Then
\begin{equation}\label{eqn:2.8}
(\EE,\FF)=(\tfrac12 \D, W^{1,2}(D))
\end{equation}
is a Dirichlet form on $L^2(D).$

\medskip

Let $\DD$ denote the class of domains in $\RR^d$ so that
the Dirichlet form $(\EE, \FF)$ of \eqref{eqn:2.8}
is regular on $L^2(\overline D)$; that is, $W^{1,2}(D) \cap C_c (\overline D)$
is both dense in $(W^{1, 2}(D), \| \cdot \|_{1,2})$ and in
 $(C_\infty (\overline D), \| \cdot \|_\infty )$.
 Here $C_\infty (\overline D)$ is the space of continuous functions on $\overline D$
 that vanishes at infinity.
 A domain $D\subset \RR^d$ is in $\DD$
  if one of the following two conditions
hold:
\begin{description}
\item{(i)} $D$ has continuous boundary (see Theorem 2 on p.14 of \cite{Maz});
that is, for every $z\in \partial D$, there is some $r>0$ so that
 $B(z, r) \cap D$ is the domain lying above the graph of a continuous function;

 \item{(ii)} $D$ is an {\it extendable domain} in the sense that
 $W^{1,2}(D)=W^{1,2}(\RR^d)\big|_D$; in other words, every  function
$u\in W^{1,2}(D)$ admits a function $\wh u\in W^{1,2}(\RR^d)$ such
that  $\wh u=u$ a.e. on $D$. When $D$ is an extendable domain,
$C_c^\infty(\overline{D})=C_c^\infty(\RR^d)\big|_{\overline{D}}$ is
a core for $(\EE, \FF)$.  Here $C_c^\infty (\overline D)$ is the
space of infinitely differentiable functions with compact support on
$\overline D$.
 \end{description}

\begin{defn} \label{D:2.1}
{\rm A domain $D$ is called a {\it locally  uniform domain}  if there are
$\delta\in (0, \infty ] $ and $C >0$ such that for every $x, y \in
D$ with $|x-y|<\delta$, there is a rectifiable curve $\gamma$ in $D$ connecting $x$ and
$y$ with $\hbox{length}(\gamma) \leq C |x-y|$ and moreover
$$\min \{ |x-z|, \, |z-y| \} \leq
 C \hbox{dist}(z, D^c) \qquad \hbox{for every } z\in \gamma.
 $$
 A domain is said to be a {\it uniform domain} if the above property holds with
 $\delta =\infty$.
  }
 \end{defn}

\medskip

  The above definition is taken from V\"ais\"al\"a \cite{Va}, where
 various equivalent definitions are discussed.
 Uniform domain and locally uniform domain are also called $(\eps, \infty )$-domain
 and $(\eps, \delta)$-domain, respectively,
 in \cite{Jo}.
 For example, the classical van Koch snowflake
domain in the conformal mapping theory is a uniform domain in
$\RR^2$.  Note that every bounded Lipschitz domain is uniform, and
every {\it non-tangentially accessible domain} defined by Jerison
and Kenig in \cite{JK} is a  uniform domain (see (3.4) of
\cite{JK}),
 while every Lipschitz domain is an $(\eps, \delta)$-domain.
 Let $S^{d-1}$ denote the unit sphere $\{x\in \RR^d: |x|=1\}$ in $\RR^d$
 and $(r, \theta)$ the polar coordinates in $\RR^d$.
 It is easy to check that a truncated infinite cone $C_{A, a}:=\{(r, \theta): \ r>a \hbox{ and } \theta \in A\}$ in $\RR^d$    for any connected open set $A\subset
S^{d-1}$ with Lipschitz boundary is a uniform domain.
However, the boundary of a uniform domain can be highly
nonrectifiable and, in general, no regularity of its boundary can be
inferred (besides the easy fact that the Hausdorff dimension of the
boundary is strictly less than $n$).
 It is known that every
 locally uniform domain
$D$ admits a bounded linear operator
\begin{equation}\label{e:extension}
 T:  \left( W^{1,2}(D), \, \| \cdot \|_{1,2} \right) \to \left(
 W^{1,2}(\RR^d), \,  \| \cdot \|_{1,2} \right)
\end{equation}
 (see   \cite{Jo}) and so in particular it is an extendable domain.

 \medskip

When $D\in \DD$, the Dirichlet form \eqref{eqn:2.8} is regular and strongly local
on $L^2 (\overline D)$. Its associated diffusion process is called symmetric
(or normally) reflected
Brownian motion (RBM) on $\overline D$. When $D=\RR^d$, then the
associated process is the $d$-dimensional standard Brownian motion
and hence (\ref{eqn:2.8}) is transient if $n\ge 3$ and recurrent if
$d=1,2.$ By the following lemma, (\ref{eqn:2.8}) is recurrent for any
$D\in \DD$ whenever $d=1,2.$

\begin{lemma}\label{lem:2.1}\
Suppose $D\in \DD$ and let $X=(X_t,\P_x)$ be RBM on $\overline D$.
\begin{description}
\item{\rm (i)}  RBM $X$ on $\overline D$
 is conservative in the sense that $\P_x(\zeta=\infty)=1$
  for q.e. $x\in \overline D$,
  where $\zeta$ is the life time of $X$.
   Moreover
   $X$ (or equivalently, the Dirichlet form $(\EE, \FF)$ of \eqref{eqn:2.8})
   is $m$-irreducible.

\item{\rm(ii)}\ Suppose $D_1,\;D_2\in \DD$ and $D_1\subset D_2.$
Let $X^{(i)}$ be RBM on $\overline D_i$, $i=1, 2$.
If $X^{(1)}$ is transient, then so is $X^{(2)}$. If $X^{(2)}$ is recurrent,
then so is $X^{(1)}$.
\end{description}
\end{lemma}

\pf\ (i)\
The conservativeness of $X$ follows from
\cite[Example 5.7.1]{FOT}.  The transition function of $X$ dominates
the transition function of the absorbed Brownian motion on $D$ which
is obtained from the $n$-dimensional standard Brownian motion by
killing upon leaving $D$ and is known to have a strictly positive
transition density on $D.$  Therefore $X$ is irreducible.

\noindent
{\rm(ii)}\ If we assume the recurrence of (\ref{eqn:2.8})
for $D_2$, there exist $u_n\in W^{1,2}(D_2),\ n\ge 1,$ which are
$\D$-Cauchy and convergent pointwise to $1$ on $D_2.$  Then
$v_n=u_n\big|_{D_1},\ n\ge 1,$ satisfy the same properties for
$D_1,$ yielding the recurrence of (\ref{eqn:2.8}) for $D_1.$  This
implies the first assertion of (ii) because of the irreducibility
proven in (i). \qed

Let us denote by
$(W^{1,2}_e(D), \EE)$
  the extended Dirichlet space of
(\ref{eqn:2.8}). We may call $W^{1,2}_e(D)$ the {\it extended
Sobolev space} of order
 $(1, 2)$
on $D.$
We denote by $\HH(D)$ the space of all harmonic functions on $D$ with finite Dirichlet integral.

\begin{lemma}\label{lem:2.11}\ $\dis W_e^{1,2}(D)\subset \bl(D)$ and
$\EE(u,u)=\frac12 \D(u,u)\ \hbox{\rm for}\ u\in W_e^{1,2}(D).$
If we denote by $\HH^*(D)$ the collection of functions
in $\bl(D)$ that is $\D$-orthogonal to all functions in $W^{1,2}_e(D)$, then
\begin{equation}\label{eqn:2.10}
\HH^*(D)\subset \HH(D).
\end{equation}
\end{lemma}

\pf\ For $u\in W^{1,2}_e(D),$ there is a sequence $\{u_n\}\subset
W^{1,2}(D)$ which is ${\bf D}$-Cauchy and convergent to $u$ a.e.
${\bf D}(u_n,u_n)$ then converges to $\EE(u,u).$  By $(\bl.1),$
there exist $v\in \bl(D)$ and constants $c_n$ such that $\{u_n\}$ is
$\bf{D}$-convergent to $v$ and the sequence $\{u_n+c_n\}$ is
convergent to $v$ in $L_{loc}^2(D).$  By choosing a subsequence if
necessary, we may assume that the latter sequence converges to $v$
a.e.  Then $\lim_{n\to\infty}c_n=c$ exists, $u=v-c$ and
consequently, $u\in \bl(D)$ and $\EE(u,u)=\frac12{\bf D}(u,u).$

If $u\in {\bl}(D)$ is ${\bf D}$-orthogonal to all functions in $W^{1,2}_e(D),$ then, since $C_c^\infty(D)\subset W^{1,2}(D),$ we have
\[(\Delta u,f)=-{\bf D}(u,f)=0,\quad \forall f\in C_c^\infty(D),\]
which implies that $\Delta u=0,$ namely, (a version of) $u$ is harmonic on $D.$ \qed

\medskip
Informally, $\HH^*(D)$ is the space of harmonic functions on $D$
having finite Dirichlet energy with zero normal derivative on the
boundary $\partial D$.

\medskip

If (\ref{eqn:2.8}) is transient, then we see from
(\ref{eqn:2.1}) that the extended Sobolev space $(W^{1,2}_e(D), \frac12\D)$ is a real Hilbert space
and in particular,
\begin{equation}\label{eqn:2.9}
u\in W^{1,2}_e(D) \ \hbox{ having } \ \D(u,u)=0\ \Longrightarrow\
u=0\ {\rm a.e.},
\end{equation}
which means
that the Hilbert space $(W^{1,2}_e(D),\frac12\D)$ is isometrically imbedded
into a closed subspace of
$(\dot{\bl}(D),\frac12{\bf D})$ by the canonical map $\bl(D)\mapsto
\dot{\bl}(D).$  Accordingly, we have

\begin{prop}\label{prop:2.1}
\ Assume that the Dirichlet form {\rm(\ref{eqn:2.8})} on $L^2(D)$ is
transient.  Then
$(W^{1,2}_e(D),\frac12\D)$ can be regarded as a closed linear subspace of the Hilbert space
$(\dot\bl(D),\frac12 \D)$ by the canonical map $\bl(D)\mapsto
\dot\bl(D).$  \\
$\bl(D)$ is a linear space spanned by $W_e^{1,2}(D)$ and
$\HH^*(D)$.
\end{prop}

The next lemma is well-known (cf. \cite{B}).

\begin{lemma}\label{lem:2.12}\ When $D=\RR^d,$ $\HH^*(\RR^d)$ consists of all constant functions.
  If $d\ge 3,$ then
\begin{equation}\label{eqn:2.101}
\bl(\RR^d)\ \hbox{\rm is a linear space spanned by}\ W^{1,2}_e(\RR^d)\ {\rm
and}\ \hbox{\rm constant functions}.
\end{equation}
\end{lemma}

See \cite[Example 1.5.3]{FOT} for a proof, where however an
incorrect statement that \lq the space $W^{1,2}_e(\RR^d)$ is
obtained from $\bl(\RR^d)$ by removing non-zero constant functions
\rq was made.  It should be corrected in the manner of
(\ref{eqn:2.101}).
See also Remark \ref{rem:2.1} below.
\medskip

 We shall be concerned with a finite measure $m(dx)=m(x)dx$
on $D$ with a density function $m(x)$ satisfying

\medskip
\noindent {\bf (A.1)}\ $\dis m(x)>0$ for every
 $ x\in \overline D$ and $m\in C_b(\overline D) \cap L^1(D)$.

\medskip
We then consider the form defined by
\begin{equation}\label{eqn:2.11}
(\EE^*,\FF^*)=\left(\frac12{\bf D},\bl(D)\cap L^2(D;m)\right),
\end{equation}
which is obtained just by replacing $L^2(D)$ with $L^2(D;m)$ in (\ref{eqn:2.7}) and (\ref{eqn:2.8}).

\begin{prop}\label{prop:2.2}\
 The symmetric form $(\EE^*, \FF^* )$ of \eqref{eqn:2.11}
is a recurrent Dirichlet form on $L^2(D;m).$  Its extended Dirichlet space
$(\EE^*, \FF_e^* )$ coincides with the space $( \frac12{\bf D}, \, \bl(D))$.
\end{prop}

\pf\ Since the convergence in $L^2(D;m)$ implies the convergence in
$L_{loc}^2(D),$ (\ref{eqn:2.11}) can be readily seen to be a closed
symmetric form on $L^2(D;m).$  Its Markovian property is an
immediate consequence of $(\bl.2)$.  It satisfies the recurrence
condition (\ref{eqn:2.3}) because $m$ is a finite measure on $D.$ We
further have
\begin{equation}\label{eqn:2.12}
u\in \FF_e^* \ \hbox{ having } \  \EE^*(u,u)=0\ \Longrightarrow u\
\hbox{\rm is constant}\ a.e.
\end{equation}
To see this, suppose $u\in \FF_e^*$ with $ \EE^*(u,u)=0$ and put
$u_\ell=\varphi^\ell\circ u$ with the normal contraction
\begin{equation}\label{e:varphi}
 \varphi^\ell(t):=((-\ell)\vee t)\wedge \ell, \qquad \ell\in \N.
\end{equation}
 Then
$u_\ell\in L^2(E;m)\cap\FF_e^*=\FF^*$ and hence
$\EE^*(u_\ell,u_\ell)=\frac12\D(u_\ell,u_\ell)=0.$  Therefore
$u_\ell$ is a constant and we get (\ref{eqn:2.12}) by letting
$\ell\to\infty.$

Denote by $\dot{\FF}_e^*$ the quotient space of $\FF_e^*$ by the
subspace of constant functions. Just as in the proof of the
preceding proposition but using
\eqref{eqn:2.12}
  in place of (\ref{eqn:2.9}), we conclude that the
space $(\dot{\FF}_e^*,\EE^*)$ is isometrically embedded into the
space $(\dot{\bl}(D),\tfrac12{\bf D}).$

Take any $u\in \bl(D)$ and put $u_\ell=\varphi^\ell\circ u$ as above.
By $(\bl.2),$ $u_\ell\in \bl(D)$ and
\begin{equation}\label{e:diri}
   {\bf D}(u-u_\ell,u-u_\ell)=\int_{\{x:|u(x)|>\ell\}}|\nabla u|^2dx\ \to\ 0,\quad \ell
   \to \infty.
\end{equation}
Since $u_\ell\in \FF^*$ and $u_\ell$ converges to $u$ pointwise, $u$
must be an element of $\FF_e^*.$  Hence the above isometric
embedding is an onto map and so $\FF_e^*=\bl(D).$
\qed

When the Lebesgue measure of the domain $D$ is finite, then we can
take
 the function $m\equiv 1$
in (\ref{eqn:2.11}) in reducing $\FF^*$ to $W^{1,2}(D).$  Hence

\begin{corollary}\label{cor:2.1}\ If the domain $D$ is of finite Lebesgue measure, then
 \hfill \break
$(\tfrac12 \D, W^{1,2}(D))$  is a recurrent Dirichlet form on $L^2(D)$ and
\begin{equation}\label{eqn:2.13}
W^{1,2}_e(D)=\bl(D).
\end{equation}
\end{corollary}

We will show in Theorem \ref{th:2.1} below that $\bl (D)$ is the reflected Dirichlet
space of $(\tfrac12 \D, W^{1,2}(D))$ for any $D\in \DD$ in any dimension $d\geq 1$.

\medskip

In what follows, we are concerned with the following condition on the domain $D\subset \RR^d$:

\medskip
\noindent
{\bf(A.2)}\ $D\in \DD$ and the Dirichlet form (\ref{eqn:2.8}) on $L^2(D)$ is transient.

\medskip

This condition forces dimension $d$ to be greater than or equal to $3$ and $D$
 to
 have infinite Lebesgue measure
 in view of Lemma \ref{lem:2.1} and Corollary \ref{cor:2.1}.
We shall see from Theorem \ref{T:2.7} below that,
if $d\ge 3$ and $D\in \DD$ contains
 an unbounded uniform domain,
then $D$ satisfies the condition {\bf (A.2)}; in other words,
RBM on $\overline D$ is transient.
This result is almost sharp since
RBM $X$ on an infinite cylinder
$D=\{x=(x_1,x_2, \cdot, x_d)\in \RR^d:  \, \sum_{i=2}^{d} x_k^2 < 1\}$ where $d\geq 3$
is recurrent because $X$ is  direct product of the one dimensional Brownian motion and
 the reflecting Brownian motion on the closed unit ball.
Note that an infinite cylinder is {\it not} a uniform domain.

We first prepare a proposition.

\begin{prop}\label{P:2.6} Let $D$ be an unbounded uniform domain in $\RR^d$ with $d\geq 3$.
and $f\in W^{1,2}_e(\RR^d)$. If $f=c$ a.e. on $D$ for some constant $c$, then $c=0$.
\end{prop}

\pf Without loss generality, we assume that the origin $0\in D$. Since
$D$ is unbounded, there is a sequence $\{x_k, k\geq 0\}$ so that $|x_k|=2Ck$,
where $C\geq 1$ is the constant in the Definition \ref{D:2.1} for the uniform domain $D$.
By definition of uniform domain, for every $k\geq 1$,
there is a rectifiable curve $\gamma_k$ connecting $0$ and $x_k$
with $2Ck=|x_k|\leq\ {\rm length}(\gamma_k)\leq C|x_k|=2C^2 k$ and
$$
\min\{|z|, |z-x_k|\}
 \leq C \, \hbox{dist} (z, D^c) \qquad \hbox{for every }
z\in \gamma_k.$$
 Choose $y_k\in \gamma_k$ such that $|y_k|=|y_k-x_k|.$
 Then clearly $|y_k|\geq |x_k|/2=Ck,\ 2|y_k|\le\ {\rm length}(\gamma_k)\le 2C^2 k$ and further
$$ \hbox{dist} (y_k, D^c) \geq (1/C) (C k)=k.
$$
Accordingly $D$ contains a sequence of open balls $\{B(y_k, k), \, k\geq 1\}$ with
$C k\leq |y_k|\leq C^2k$.

Let $X=(X_t,\P_x)$ be the Brownian motion in $\RR^d$, which is transient as $d\geq 3$,
and $e_1:=(1, 0, \cdots, 0)\in \RR^d$.
Observe that by the Brownian scaling and the rotation invariance of the Brownian motion $X$,
\begin{eqnarray*}
 && \P_0 (X_t\in D\ \hbox{for some}\ t\ge Ck)\\
 & \geq & \P_0 ( X_t\in B(y_k, k) \ \hbox{for some}\ t\ge Ck)\\
&\ge& \P_0 ( X_t\in B(y_k, k) \ \hbox{for some}\ t\ge |y_k|)\\
&=& \P_0 ( X_t\in B (e_1, k/|y_k|)\ \hbox{for some}\ t\ge 1/|y_k| ) \\
&\geq & \P_0 ( X_t\in B (e_1, C^{-2})\ \hbox{for some}\ t\ge C^{-1}) =:p_0>0.
\end{eqnarray*}
Consequently
 \begin{equation}\label{e:hit}
 \P_0 (
  \cap_{n\geq 1} \cup_{k\geq n}
 \{X_t\in D\ \hbox{for some}\ t\ge Ck\})\ge p_0.
\end{equation}

Let $\partial$ be the one-point compactification of $\RR^d$.  As is well-known,
\begin{equation}\label{eqn:escape}
\P_x\left( \lim_{t\to\infty}X_t=\partial\right) = 1,\quad x\in \RR^d.
\end{equation}
Any $f \in W^{1,2}_e(\RR^d)$ admits a version $\wt f$ which is quasi continuous in the restricted sense with respect to the
transient extended Dirichlet space \\
$(W_e^{1,2}(\RR^d),\frac12{\bf D})$, namely,
there exists a  decreasing sequence of
open subsets $\{G_n\}$ of $\RR^d$ such that ${\rm Cap}_{(0)}(G_n)\to 0,\ n\to \infty,$ and the restriction of $\wt f$ to each set $(\RR^d \setminus G_n)\cup \{\partial\}$ is continuous there if we set $\wt f(\partial)=0.$  By (2.9) of \cite{CF1}, it then holds that $\P_x(\sigma_{G_n}=\infty\ \hbox{for some}\ n\ge 1)=1,\ {\rm q.e.}\ x\in \RR^d,$ which combined with \eqref{eqn:escape}
yields that $\varphi(x)=1$ for q.e. $x\in \RR^d$ where $\varphi(x)=\P_x\left(\lim_{t\to \infty} \wt f(X_t) = 0\right).$  Hence, by the Markov property,
\begin{equation}\label{eqn:zero}
\P_0\left(\lim_{t\to\infty}\wt f(X_t)=0\right)=(2\pi)^{-d/2}\int_{\RR^d} e^{-|x|^2/2}\varphi(x)dx =1.
\end{equation}

Suppose $f\in W_e^{1,2}(\RR^d)$ equals a constant $c$ a.e. on $D$, then $\wt f=c$ q.e. on $D,$ and we can deduce from \eqref{e:hit} and \eqref{eqn:zero} that $c=0$. \qed

\begin{remark}\label{rem:2.1}\rm\ We can make the statement \eqref{eqn:2.101} more specific as follows: When $d\ge 3,$ $W_e^{1,2}(\RR^d)$ is obtained from the quotient space $\dot \bl(\RR^d)$ by choosing from each equivalence class a special function whose quasi continuous version has zero limit at $\partial$ in the sense of \eqref{eqn:zero}.
\end{remark}

\begin{thm}\label{T:2.7} Assume that $d\ge 3$ and $D$ is a domain in $\DD$ that contains
an unbounded uniform domain. Then RBM on $\overline D$ is transient; in other words,
 $D$ satisfies condition {\bf (A.2)}.
\end{thm}

\pf In view of the comparison result Lemma \ref{lem:2.1}(ii),
we may assume, without loss of generality, that $D$ is an unbounded uniform domain.

We first show that $W^{1,2}_e(D)=W^{1,2}_e(\RR^d)\big|_D$.
By \cite[Theorem 4.13]{HK},
  there is a linear bounded extension
operator $S : \bl (D)\to \bl (\RR^d)$ in the following sense: for every
$u\in \bl (D)$, $S u\in \bl (\RR^d)$ with $Su=u$ a.e. on $D$, moreover
\begin{equation}\label{e:ext2}
 \| \nabla (Su)\|_{L^2(\RR^d)} \leq M \, \| \nabla u\|_{L^2(D)}
\qquad \hbox{for every } u \in \bl (D).
\end{equation}
For $u\in W^{1,2}(D)$,  $Su\in \bl (\RR^d)$ and so by Lemma \ref{lem:2.12},
$Su=v_0+c$ for some $v_0 \in W^{1,2}_e(\RR^d)$ and constant $c$.
On the other hand, by \eqref{e:extension}, there is
 $Tu\in W^{1,2}(\RR^d)$ so that $Tu =u $ a.e. on $D$.
It follows that $Tu-v_0 \in W^{1,2}_e(\RR^d)$ and $Tu-v_0 = c$ on $D$.
Proposition \ref{P:2.6} now implies that $c=0$.
This proves that the extension
operator $S$ maps  $W^{1,2}(D)$ into $W^{1,2}_e(\RR^d)$.
Now for $u\in W^{1,2}_e(D)$, there is a $\D$-Cauchy sequence $\{u_k, k\geq 1\}\subset
 W^{1,2}(D)$ so that $u_k\to u$ a.e. on $D$. By \eqref{e:ext2} and what we just established,
 $\{Su_k, k\geq 1\}$ is a $\D$-Cauchy sequence in $W^{1,2}_e(\RR^d)$.
 Since $(W^{1,2}_e (\RR^d), \D)$ is a Hilbert space, $Su_k$ is $\D$-convergent to some
 $f\in W^{1,2}_e(\RR^d)$. Taking a subsequence if necessary, $Su_k$ converges to $f$ a.e. on
 $\RR^d$. It follows then $u=f$ a.e. on $D$. This proves that
 $W^{1,2}_e(D)=W^{1,2}_e(\RR^d)|_D$, as  clearly $W^{1,2}_e(\RR^d)|_D\subset W^{1,2}_e(D)$.
(The above also implies that $S$ maps $W^{1,2}_e(D)$ into $W^{1,2}_e(\RR^d)$.
 This is because for $u\in W^{1,2}_e(D)$, $Su \in
  \bl (\RR^d)$ and so $Su=v_0+c$
 for some $v_0\in W^{1,2}_e(\RR^d)$ and constant $c$. On the other hand, there is
 some $f\in W^{1,2}_e(\RR^d)$ so that $f=u$ a.e. on $D$. Thus we have
 $f-v_0=c$ on $D$. Since $f-v_0\in W^{1,2}_e(\RR^d)$, we have by Proposition \ref{P:2.6}
 that $c=0$. This proves that $Su=v_0\in W^{1,2}_e(\RR^d)$.)

 To show that RBM on $\overline D$ is transient,
 by \eqref{eqn:2.2} and Lemma \ref{lem:2.11},
 it suffices to show that if $u\in W^{1,2}_e(D)$ has
 $\D (u, u)=0$, then $u=0$ on $D$.
 Suppose $u\in W^{1,2}_e(D)$ and $\D (u, u)=0$.
Clearly, $u=c$ a.e. on $D$  for some constant $c$.
We know from the above that there is
 $u_0\in W^{1,2}_e(\RR^d)$ so that $u_0=u=c$ a.e. on $D$.
It follows from Proposition \ref{P:2.6} that $c=0$ and so $u=0$. \qed

\section{Unique symmetric extension and the space of BL-functions}

We start this section with a domain $D\in \DD$, together with a measure
$m(dx)=m(x)dx$ on $D$ whose density $m(x)$ is strictly positive for every $x\in \overline D$
and satisfies $m\in C_b (\overline D)$.
Consider the following symmetric form
\begin{equation}\label{eqn:2.14}
(\EE^{(0)},\FF^{(0)})=\left(\frac12{\bf D}, \, W^{1,2}_e(D)\cap
L^2(D;m)\right),
\end{equation}
obtained by just replacing $\bl(D)$ with $W^{1,2}_e(D)$ in
(\ref{eqn:2.11}).
As for (\ref{eqn:2.11}), the above form can be readily
checked to be a Dirichlet form on $L^2(D;m).$   Furthermore it is a
 regular
and strongly local
Dirichlet form on $L^2(\overline{D};m).$  The associated
diffusion
 process $Y=(Y_t,\P_x)$ on $\overline{D}$ is obtained
from the reflecting Brownian motion $X$ by a time change with
respect to the additive functional $\dis A_t=\int_0^tm(X_s)ds$ in
view of \cite[\S 6.2]{FOT}.

For a general regular
Dirichlet form $(\EE,\FF),$ the notions of
its reflected Dirichlet space $(\FF^{\rm ref}, \EE^{\rm ref})$ and
its active reflected Dirichlet form $(\EE^{\rm ref},\FF^{\rm ref}_a)$
were originally introduced by M. L. Silverstein (see \cite{S1}) and
have been further studied in
\cite{C} when $(\EE, \FF)$ is transient.
These notions are well
defined for the present Dirichlet form (\ref{eqn:2.14}).

Assume that $m\in C_b( \overline D)$ with $m(x)>0$ for every $x\in \overline D$
and let $(\EE^{(0)}, \FF^{(0)})$ be the Dirichlet form on $L^2(D; m)$ defined by
\eqref{eqn:2.14}.
Following \cite{C}, let
\begin{eqnarray*}
 \FF^{\rm ref} &:=&\left\{ f : \ f_k:= \varphi^k\circ f \in \FF^{(0)}_{\rm loc}
 \hbox{ and }
 \sup_{k\geq 1} {\bf D} (f_k, f_k)<\infty  \right\},\\
\FF^{\rm ref}_a &:=& \FF^{\rm ref} \cap L^2 (D; m),\\
\EE^{\rm ref} (f, f) &:=& \frac12 \int_{D} |\nabla f(x)|^2 dx = \frac12
 {\bf D}(f, f) \qquad \hbox{for } f\in \FF^{\rm ref} .
  \end{eqnarray*}
We call  $(\EE^{\rm ref}, \FF^{\rm ref})$ and $(\EE^{\rm ref}, \FF^{\rm ref}_a)$
the reflected Dirichlet form and the active
reflected Dirichlet form, respectively, of $(\EE^{(0)},\FF^{(0)})$
 on $L^2(D; m)$.
In \cite{C}, the above notions are defined for transient $(\EE^{(0}, \FF^{(0)})$
and it is proved there that $(\EE^\rf, \FF^\rf_a)$ is a Dirichlet form on $L^2(D; m)$.
But these notions can also be defined when  $(\EE^{(0}, \FF^{(0)})$ is recurrent.

\begin{thm}\label{th:2.1}\
Let $D\in \DD$ be a domain in $\RR^d$ with $d\geq 1$.
Assume that $m\in C_b( \overline D)$ with $m(x)>0$ for every $x\in \overline D$.
Then the Dirichlet form $(\EE^*,\FF^*)$ defined by
 \eqref{eqn:2.11}
  with $m(dx)=m(x)dx$ is
 the active reflected Dirichlet form of the
Dirichlet form $(\EE^{(0)},\FF^{(0)})$ given by
 \eqref{eqn:2.14}
 and $\bl (D)$ is the reflected Dirichlet space of $(\EE^{(0)},\FF^{(0)})$.
\end{thm}

\pf
For $f\in \FF^\rf\cap L^\infty (D)$,   $f\in \FF^{(0)}_{loc}$ and so
$f\in \bl (D)$. For general $f\in \FF^\rf$ and $k\geq 1$,
 define $f_k=\varphi^k\circ f$ for $k\geq 1$,
where $\varphi^k$ is defined in \eqref{e:varphi}.
From above we have $f_k\in \bl (D)$ with $\EE^\rf (f_k, f_k)= \frac12 \D (f_k, f_k)$.
Since $\sup_{k\geq 1} \D (f_k, f_k)<\infty$, there is a Ces\`aro mean $\{g_k, k\geq 1\}$
of $\{f_k, k\geq 1\}$ that is $\D$-Cauchy. Observe that $g_k\in \FF^\rf\cap L^\infty (D)
\subset \bl (D)$. By (BL.1), there is $u\in \bl (D)$ and a sequence of constants
$\{c_k, k\geq 1\}$ so that $g_k$ is $\D$-convergent to $u$ and $g_k+c_k$ converges to
$u$ locally in $L^2(D)$ as $k\to \infty$. But $\lim_{k\to\infty} g_k=f$ on $D$.
This implies that  $\lim_{k\to \infty} c_k=c$ and  $f+c =u$.
Therefore $f\in \bl (D)$ and so $\FF^\rf\subset \bl (D)$.

Now for $f \in \bl (D)$, clearly
$f_k=\varphi^k\circ f \in \bl (D) \cap L^2_{\rm loc}(D)\subset \FF^{(0)}_{\rm loc}$
 for every $k\geq 1$. Since ${\bf D} (f_k, f_k) \leq {\bf D}(f, f)$,
 we conclude that $f\in \FF^{\rm ref}_a$. This completes the proof that
 $ \bl (D)   =\FF^{\rm ref}  $. It then follows immediately that
 $\FF^*=\bl (D)\cap L^2(D; m)= \FF^\rf\cap L^2(D; m)=\FF^\rf_a$. \qed

 Taking $m(x)\equiv 1$ on $D$ in Theorem \ref{th:2.1} yields in particular that
$\bl (D)$ is the reflected Dirichlet space for $(\tfrac12 \D, W^{1,2}(D) )$ for every
$D\in \DD$ in any dimension.

\medskip

We shall next consider a domain $D\in \DD$
     satisfying the transience condition {\bf(A.2)}  and the function $m$
satisfying condition {\bf(A.1)}.
Then $d\ge 3$ and $D$ is of infinite Lebesgue measure as is noted in \S 2.
 Abusing the notation a bit, we also use $\partial$ to denote
 the point at infinity of $\overline{D}$,  namely,
  we let $\overline{D}^*=\overline{D}\cup \{\partial\}$
   be the one point compactification of $\overline{D}.$
Under condition {\bf(A.2)}, the reflecting Brownian motion $X=(X_t, \P_x)$ on $\overline{D}$ satisfies
\begin{equation}\label{eqn:wanderout}
 \P_x(\lim_{t\to\infty}X_t=\partial)=1 \qquad  \hbox{for q.e. }
  x\in \overline{D}
\end{equation}
on account of \cite[Theorem 2.4]{CF1}.
In this case, the Dirichlet form $(\EE^{(0)},\FF^{(0)})$ of
\eqref{eqn:2.14} is transient
 and its associated diffusion process, time-changed Brownian motion
 $Y=(Y_t,\P_x)$ on $\overline{D}$  has
  lifetime $\zeta$
 satisfying
$\P_m(\zeta<\infty)>0$. This is because, for the resolvent $\{R_\alpha;
\alpha>0\}$ of $Y,$ $R_\alpha 1$ is an element of $\FF^{(0)}$ but it
can not be a positive constant $1/\alpha$
 by virtue of \eqref{eqn:2.9}.
Moreover, $Y_t$ approaches to $\partial$ as $t\to \zeta$ in view of
\eqref{eqn:wanderout}.

\medskip

Let us further make the following assumption on the domain $D$:

\medskip
\noindent
{\bf(A.3)} \ $\HH^*(D)$\ \hbox{\rm consists
of
constant functions on}\ $D.$

\medskip
\noindent Here $\HH^*(D)$ is the space of harmonic functions on $D$
appearing in
Lemma \ref{lem:2.11}. Thus under the transience condition {\bf (A.2)},
 {\bf(A.3)} amounts to
assuming that
\begin{equation}\label{eqn:3.1}
\bl(D)\ \hbox{\rm is a linear space spanned by}\ W^{1,2}_e(D)\ {\rm
and}\ \hbox{\rm constant functions}.
\end{equation}

We will show in Theorem \ref{T:3.5} below
that conditions {\bf (A.2)} and {\bf (A.3)} hold for every unbounded
uniform domain.
 Moreover
Proposition \ref{prop:3.1} will imply
that the condition {\bf (A.3)} is
satisfied by any domain $D\in \DD$ that is the complement of a
compact set.

Under the stated three conditions, we shall compare the form $(\EE^*,\FF^*)$
defined by (\ref{eqn:2.11}) and
the $(\EE^{(0)},\FF^{(0)})$ defined by (\ref{eqn:2.14}).
They are Dirichlet forms on $L^2(D;m)$.
 and related to each other by
\begin{equation}\label{eqn:3.2}
\left\{
\begin{array}{ccl}
\FF^*&=& \left\{u=u_0+c:u_0\in \FF^{(0)},\ c\ \hbox{is constant} \right\}, \\
\EE^*(u,u)&=&\EE^{(0)}(u_0,u_0)=\frac12\D(u_0,u_0).
\end{array}
\right.
\end{equation}
in view of Proposition \ref{prop:2.1}.

Let $m^*$ be the extension of $m$ from $D$ to $\overline{D}^*$ obtained by setting
$$ m^*(\partial D\cup \{\partial\})=0.
$$
By identifying $L^2(D;m)$ with $L^2(\overline{D}^*;m^*)$,
 we can regard $(\EE^*,\FF^*)$ as a Dirichlet form on $L^2(\overline{D}^*;m^*).$

\begin{thm}\label{th:3.1}\
Assume that conditions {\rm\bf(A.1), (A.2)} and {\rm\bf(A.3)} hold.
\begin{description}
\item{\rm(i)}\ $(\EE^*,\FF^*)$ is a recurrent, strongly local and
regular Dirichlet form on $L^2(\overline{D}^*,m^*).$

\item{\rm(ii)}\ The associated diffusion process $Y^*=(Y^*_t,\P_x^*)$
on $\overline{D}^*$ is a conservative extension of the time changed
transient reflected Brownian motion $Y$ on $\overline{D}$ and
satisfies
\begin{equation}\label{eqn:3.3}
\P_x^*(\sigma_{\partial}<\infty)=1 \qquad \hbox{for q.e.}\ x\in
\overline{D}^*.
\end{equation}
\end{description}
\end{thm}

\pf\ (i)\ The recurrence condition (\ref{eqn:2.3}) is trivially
satisfied by $(\EE^*,\FF^*).$   If we let $\CC=\{u+c: u\in
C_c^\infty(\overline{D}),\ c\ \hbox{\rm is constant}\},$ then
$\CC\subset C(\overline{D}^*)$ and $\CC$ is readily seen to be a
core of the Dirichlet form $(\EE^*,\FF^*).$  The strong locality of
$(\EE^*,\FF^*)$ can be proved in the same way as in the proof of
Theorem 3.2 of \cite{FT}, where a Dirichlet form quite similar to
(\ref{eqn:3.2}) was studied in a rather general context.

\noindent (ii)\ Since $Y$ is associated with $(\EE^0,\FF^0),$ $Y^*$
is an extension of $Y$ from $\overline{D}$  to
$\overline{D}^*=\overline{D}\cup \{\partial\}$. The point $\partial$ is
not $m^*$-polar for $X^*$ because
$\P^*_{m^*}(\sigma_{\partial}<\infty)$ dominates the quantity
$\P_m(\zeta<\infty)>0$ for $Y.$ On the other hand, the
irreducibility of $(\EE^*,\FF^*)$ can be verified in a similar
manner to \cite{FT}. Hence (\ref{eqn:3.3})
holds
 in view of \cite[Theorem 4.6.6]{FOT}. \qed

\begin{remark}\label{rem:3.1}\rm The extended diffusion $Y^*$ of $Y$ in the above theorem can be also constructed stochastically by a darning of $Y$ at $\partial$ as in \cite{FT} and \cite{CF1}, namely, by piecing together the excursions of $Y$ around $\partial$
according to an excursion-valued Poisson point process $\{\p_t\}.$  The characteristic measure $\n$ of $\{\p_t\}$ is described as \cite[(4.4)]{FT} in terms of the transition function $\{q_t\}$ of $Y$ and the $\{q_t\}$-entrance law $\{\mu_t\}$ on $E$ uniquely determined by the equation
\[m=\int_0^\infty \mu_t dt.\]
\end{remark}

Using Theorem \ref{th:2.1}, we can further show that $Y^*$ is the
unique $m$-symmetric continuous genuine extension of $Y$ in the
following sense.

\begin{thm}\label{T:3.2}
Suppose that there is a Luzin space $E$ so that $\overline D$ is
embedded continuously in $E$ as a dense open subset and there is an
$m$-symmetric diffusion $Z$ on $E$ that admits no killings inside
and that $Y$ is the proper subprocess of $Z$ killed upon leaving
$\overline D$. Here the measure $m$ is extended to $E$ by setting
$m(E\setminus \overline D)=0$. Then the symmetric Dirichlet form of
$Z$ on $L^2(E; m)$ coincides with that of $X^*$ on $L^2(\overline D;
m)$; in other words, $Z$ under $\P_m$ has the same
finite-dimensional distributions as that of $Y^*$ under $\P_m^*$.
\end{thm}

\medskip

 Note that although $Z$ and $X^*$ live on different state spaces,
 since $m(E\setminus \overline D)=0$ and $m(\{\partial \})=0$,
 $$\P_m(Z_t\in E\setminus \overline D)=\P_m^*(Y^*_t = \partial )=0 \qquad
 \hbox{for every } t\geq 0
 $$
  and so we can compare
 the
finite-dimensional distributions of $Z$ under $\P_m$ with
those of $Y^*$ under $\P_m^*$.

\bigskip

\noindent{\bf Proof of Theorem \ref{T:3.2}.} Let $(\EE^Z, \FF^Z)$
denote the symmetric Dirichlet form of $Z$ on $L^2(E; m)$. It is
known (cf. \cite{MR}) that $(\EE^Z, \FF^Z)$ is quasi-regular on $E$
and so, by \cite{CMR}, it is quasi-homeomorphic to a regular
Dirichlet form on a locally compact metric space. Consequently, the
results established under the regular Dirichlet form framework in
\cite{FOT} and \cite{S1} apply to $Z$ and $(\EE^Z, \FF^Z)$.

Let $\LL$ and $\LL^0$ be the $L^2$-generators of the Dirichlet forms
$(\EE^Z, \FF^Z)$ on $L^2(E;m)$ and $(\EE^{(0)}, \FF^{(0)})$ on
$L^2(\overline D; m)$, respectively. Since $Y$ is the part process
of $Z$ killed upon leaving $\overline D$, for every $f\in {\rm
Dom}(\LL)$, there is a unique $\EE^Z$-orthogonal decomposition
$f=f_0+h$, where $f_0\in \FF^{(0)}$ and $h\in \FF^Z$ so that
$\EE^Z(g, h)=0$ for every $g\in \FF^{(0)}$. Note that $h$ is
harmonic with respect to the diffusion $Y$ (see \cite{FOT}). Since
$f\in {\rm Dom}(\LL)$, we have in particular
$$ (-\LL f, g)_{L^2(\overline D; m)} =\EE^Z( f, g)=\EE^{(0)}(f_0, g)
\qquad \hbox{for every } g\in \FF^{(0)}\subset \FF^Z.
$$
The above implies that $f_0\in {\rm Dom}(\LL^0)$ and $\LL^{(0)} f_0
= \LL f$. In the terminology of Silverstein \cite[{\bf 15.1} on
p.152]{S1}, $\LL$ is contained in the {\it local generator} of
$(\EE^{(0)}, \FF^{(0)})$. So by \cite[Theorem 15.2]{S1} and its
proof as well as Theorem \ref{th:2.1} above, we have
$$ \FF^Z \subset \FF^* \qquad \hbox{and} \qquad
 \EE^Z (f, f) \geq \EE^*(f, f) \  \hbox{ for } f\in \FF^Z.
 $$
Since $Z$ is a genuine extension of $Y$, $\FF^Z$ must be the same as
$\FF^*$ in view of \eqref{eqn:3.2}. For any $f=f_0+c_1\in \FF^*$
with $f_0\in \FF^{(0)}$, since $(\EE^Z, \FF^Z)$ is local and admits
no killing, we have
$$ \EE^Z(f, f)= \EE^Z(f_0, f_0)= \EE^{(0)}(f_0, f_0)= \EE^*(f, f).
$$
This proves that $(\EE^Z, \FF^Z)=(\EE^*, \FF^*)$.
 \qed

\bigskip

Finally we give sufficient conditions for {\bf (A.3)} to hold.

\begin{thm}\label{T:3.5}
Let $D\subset \RR^d$ be a uniform domain.
Then $\HH^*(D)$ consists of all constant functions on $D.$
In particular, if $d\ge 3$ and $D\subset \RR^d$ is an unbounded uniform domain, then $D$ satisfies conditions  {\rm\bf (A.2)} and {\bf (A.3)}.
\end{thm}

\pf
Again we use \cite[Theorem 4.13]{HK} which states that there is a linear bounded extension
operator $S : \bl (D)\to \bl (\RR^d)$ in the following sense: for every
$u\in \bl (D)$, $S u\in \bl (\RR^d)$ with $Su=u$ a.e. on $D$
and with bound \eqref{e:ext2}.

For $u\in \HH^*(D)\subset \bl (D)$, $Su\in \bl (\RR^d )$.
Therefore, by Lemma \ref{lem:2.12},
$Su= f_0+c$ for some $f_0\in W^{1,2}_e(\RR^d)$ and some constant $c$. Note
that $f_0|_D \in W^{1,2}_e(D)$.
Since $f_0|_D=u-c\in\HH^*(D)$, we obtain
\begin{equation}\label{e:f0}
 \D (f_0, f_0)=0,
\end{equation}
which implies that $f_0=c_1$ for some constant.  Hence $u\in \HH^*(D)$ equals a constant $c+c_1$ on $D.$

The second statement of the theorem follows from the
 the above
and Theorem \ref{T:2.7}. \qed

We can readily deduce the next proposition from the above theorem.  In particular, any unbounded domain $D\subset \RR^d,\ d\ge 3,$ in $\DD$ with compact boundary satisfies {\bf(A.3)}.

\begin{prop}\label{prop:3.1}\ If $D\subset \RR^d$, $d\geq 3$, is an unbounded domain in $\DD$ such that $D\setminus
 \overline{B(0, r)}$
is a uniform domain for some $r>0,$
then condition {\rm\bf(A.3)} holds.
\end{prop}
\pf.\ By normal contraction property, it suffices to show
that for every bounded non-negative $u\in \bl (D)$, $u=u_0+c$ for some
$u_0\in W^{1,2}_e(D)$ and a constant $c$.

We put $U=D\setminus B (0, r)$
and take $\phi \in C^\infty_c(\RR^d)$ such that $\phi =1$ on $B (0, r+1)$.
As $u\in \bl (D)\subset \bl (U)$ and $U$ is an unbounded uniform domain, we see by the above theorem that there is a constant $c$
with $u_0:=u-c\in W^{1,2}_e(U)$. We claim that $(1-\phi) u_0
\in W^{1,2}_e(U)$. This is because there is a sequence $f_n \in
W^{1,2}(U)$ such that $\sup_{n\geq 1} \int_U |\nabla f_n|^2 dx =0$
and $f_n\to u_0,\ n\to\infty,$ a.e. on $U$. Since $u_0$ is bounded,
by the normal contraction property, we may assume that $\{f_n, n\geq
1\}$ are uniformly bounded. Clearly $(1-\phi) f_n\in W^{1,2}(D)$,
$(1-\phi)f_n$ converges a.e. on $D$ to $(1-\phi) u_0$ and
 \begin{eqnarray*}
 && \sup_{n\geq 1} \int_D |\nabla ((1-\phi) f_n)|^2 dx \\
 &\leq &
2 \sup_{n\geq 1} \|f_n \|^2_{L^\infty(U)} \int_{\RR^d} |\nabla \phi|^2 dx
+ 2 \|(1-\phi) \|_\infty^2 \sup_{n\geq 1} \int_U |\nabla f_n|^2 dx <\infty.
\end{eqnarray*}
We conclude that $(1-\phi) u_0 \in W^{1,2}_e(D)$. On the other hand,
since $\phi u_0 \in W^{1,2}(D)$, we have
 $u-c=u_0 = \phi u_0 + (1-\phi ) u_0 \in W^{1,2}_e(D)$. This proves the theorem.\qed

\begin{remark}\label{rem:3.2}{\rm
The second statement of Theorem \ref{T:3.5} is ``sharp" for {\bf (A.3)} in view of the following example.

 If $D$ is an unbounded domain in $\RR^d$ with two or more
infinite branches, then condition {\bf (A.3)} may fail. Consider
$$ D=B (0, 1) \cup \Big\{x=(x_1, x_2, \cdots, x_d): \,
x_d^2 >\sum_{k=1}^{d-1}x_k^2 \Big\}
$$
in $\RR^d$ for $d\geq 3$.
Here $B(x, r)$ denote the ball centered at $x$ with radius $r>0$.
Clearly $D$ is a Lipschitz domain but it
is {\it not} a uniform domain as it has a bottleneck $B_1$.
However $D$ contains an unbounded uniform domain
$$  C_1:=B(0, 1)^c \,  \cap \, \left\{  x_d  >\sqrt{x_1^2+\cdots +x_{d-1}^2} \right\}.
$$
So by Theorem \ref{T:2.7}, condition {\bf (A.2)} holds for $D$.

We claim
that condition {\bf (A.3)} does not hold for $D$. Define
$$
 C_2 =B(0, 1)^c \, \cap \, \left\{  x_d  <-\sqrt{x_1^2+\cdots +x_{d-1}^2} \right\}.
$$
Let $f\in C^2_b (\overline D)$ be such that $f=1$ on $C_1$ and $f=2$ on $C_2$. Clearly
$f\in \bl  (D)$ and so we can write it as $f=f_0+h$ with
$f_0\in W^{1,2}_e(D)$ and $h\in \HH^*(D)$.
Let $X$ be RBM on $\overline D$, which is transient as is noted in the above.
By \eqref{eqn:wanderout}
$$ \P_x \left( \lim_{t\to \infty} f(X_t)=1 \right) \cdot  \P_x
\left( \lim_{t\to \infty} f(X_t)=2 \right )>0
 \qquad \hbox{for q.e.} x\in D.
 $$

On the other hand, $f_0$ is quasi continuous in the restricted sense with respect to the transient extended Dirichlet space $(W_e^{1,2}(D),\frac12{\bf D})$ : there exists a  decreasing sequence of
open subsets $\{G_n\}$ of $\overline D$ such that ${\rm Cap}_{(0)}(G_n)\to 0$ as
$ n\to \infty,$ and the restriction of $f_0$ to each set $(\overline D \setminus G_n)\cup
\{\partial \}$ is continuous there if we set $f_0(\partial)=0.$  By (2.9) of \cite{CF1}, it then holds that $\P_x(\sigma_{G_n}=\infty\ \hbox{for some}\ n\ge 1)=1,\ {\rm q.e.}\ x\in D,$ which combined with \eqref{eqn:wanderout} leads us to
$\lim_{t\to \infty} f_0(X_t) = 0$ $\P_x$-a.s. for q.e. $x\in D$.
Thus we have $ \lim_{t\to \infty}
h(X_t) = \lim_{t\to \infty} f(X_t)$ $\P_x$-a.s. for q.e. $x\in D$.
This yields that $h$ can not be a constant. \qed }
\end{remark}

\begin{remark}\label{rem:3.3}\rm\ In many cases where the reflecting Brownian motion on $\overline{D}$ is recurrent, we have the identification
\begin{equation}\label{e:id}
(W_e^{1,2}(D),\tfrac12\D)=(\bl(D), \tfrac12\D).
\end{equation}
In Corollary \ref{cor:2.1}, this has been verified for any domain $D\subset \RR^d,\ d\ge 1,$ with finite Lebesgue measure.  When $d\le 2,$ \eqref{e:id} holds for any domain $D\subset \RR^d.$  Here is a proof when $d=2.$

Take any function $u\in\bl(D)$ such that $|u|\le \ell$ for some constant.   Let $\psi_n\in
C_c^1(\RR_+)$ be functions satisfying
\[
\left\{
\begin{array}{lrl}
\psi_n(x)=1,&\ 0\le x<n;\quad &\psi_n(x)=0,\ x>2n+1;\\
 |\psi_n'(x)|\le \frac{1}{n},&\ n\le x\le 2n+1,\ &0\le \psi_n(x)\le 1,\ x\in \RR_+.
 \end{array}
 \right.
\]
and put $u_n(x)=u(x)\psi_n(|x|), \ x\in D.$ Then $u_n\in \bl(D)\cap
L^2(D)=
 W^{1,2}(D)$
and
\begin{eqnarray*}
&&{\bf D}(u_n,u_n)\le 2\int_D|\nabla u|^2(x)\psi_n(|x|)^2dx+2\int_Du^2(x)\psi_n'(|x|)^2dx\\
&&\le 2{\bf D}(u,u)+ 2\ell^2\int_{\{x\in \RR^2:|x|\le 2n+1\}}\psi_n'(r)^2drd\theta\\
&&\le 2{\bf D}(u,u)+\frac{2\ell^2\pi(2n+1)^2}{n^2}\le 2{\bf D}(u,u)+18\ell^2\pi.
\end{eqnarray*}
Hence a Ces\`aro mean of a subsequence of $\{u_n\}$ is ${\bf D}$-convergent.
Since $u_n$ converges to $u$ pointwise, we conclude that $u\in W_e^{1,2}(D)$ and
$\EE(u,u)=\frac12{\bf D}(u,u).$

Next take any $u\in \bl(D)$ and put $u_\ell=\varphi^\ell\circ u,\
\ell\in \mathbb{N},$ for the normal contraction $\varphi^\ell$ of
(\ref{e:varphi}).
By (\bl.2),
$u_\ell\in \bl(D)$ and we have \eqref{e:diri}, which particularly means that $\D(u_\ell, u_\ell)$ is bounded.
We have just shown that $u_\ell\in W_e^{1,2}(D)$ with
$\EE(u_\ell,u_\ell)=\frac12{\bf D}(u_\ell,u_\ell)$.

Under the general setting in the beginning of \S 2, let $(\FF_e,\EE)$ be the extended Dirichlet space of a Dirichlet form $(\EE,\FF)$ on $L^2(E;m).$  Assume that an $m$-measurable function $f$ on $E$ is finite $m$-a.e. on $E,$ $f_\ell=\varphi^\ell\circ f\in \FF_e$ for each $\ell,$ and $\sup_\ell \EE(f_\ell,f_\ell)<\infty.$  Then it can be readily shown that $f\in \FF_e.$  In the present case, we have therefore $u\in W_e^{1,2}(D)$ and $\EE(u,u)=\frac12\D(u,u).$

It is possible to show that
\eqref{e:id} holds for any $D\in \DD$ for which the RBM on $\overline D$ is recurrent.
We plan to address it elsewhere in future.
\end{remark}

\vskip 0.4truein

\noindent {\bf Zhen-Qing Chen}:

\smallskip

\noindent Department of Mathematics,
University of Washington, Seattle,
 WA 98195, USA.

\noindent Email:  {\texttt zchen@math.washington.edu}

\medskip
\noindent {\bf Masatoshi Fukushima}:

\smallskip

\noindent
 Branch of Mathematical Science,
Osaka University, Toyonaka, Osaka 560-0043, Japan.

\noindent Email: {\texttt fuku2@mx5.canvas.ne.jp}

\end{document}